\documentclass[twoside,12pt]{article}
\usepackage{amssymb,amsmath,bm,euscript}
\usepackage{graphicx,color,caption}

\usepackage[colorlinks=true]{hyperref}

\usepackage{algorithm,algpseudocode}
\usepackage{braket}

\newtheorem{proposition}{Proposition}[section]
\newtheorem{theorem}[proposition]{Theorem}

\newcommand{\qed}{\hphantom{.}\hfill $\Box$\medbreak}

\def\F{\mathcal{F}}

\def\e{{\bf e}}
\def\x{{\bf x}}
\def\y{{\bf y}}
\def\z{{\bf z}}

\def\uu{{\bf u}}
\def\vv{{\bf v}}
\def\0{{\bf 0}}

\title{\bf{Nilpotent Decomposition of Solvable Lie Algebras}\thanks{This research was supported by the Hong Kong
    Research Grant Council (Grant No.  PolyU 15300715, 15301716 and 15300717). }}
\author{ \hspace{1mm} Liqun Qi\thanks{Department of Applied
    Mathematics, The Hong Kong Polytechnic University, Hung Hom,
    Kowloon, Hong Kong; ({\tt liqun.qi@polyu.edu.hk}).}
    }

\begin{document}
\date{\today}
\maketitle

\begin{abstract}  Semisimple Lie algebras have been completely classified by Cartan and Killing.  The Levi theorem states that every finite dimensional Lie algebra is isomorphic to a semidirect sum of its largest solvable ideal and a semisimple Lie algebra.   These focus the classification of solvable Lie algebras as one of the  main challenges of Lie algebra research.   One approach towards this task is to take a class of nilpotent Lie algebras and construct all extensions of these algebras to solvable ones.  In this paper, we propose another approach, i.e., to decompose a solvable nonnilpotent Lie algebra to two nilpotent Lie algebras which are called the left and right nilpotent algebras of the solvable algebra.   The right nilpotent algebra is the smallest ideal of the lower central series of the solvable algebra, while the left nilpotent algebra is the factor algebra of the solvable algebra and its right nilpotent algebra.   We show that the solvable algebras are decomposable if its left nilpotent algebra is an Abelian algebra of dimension higher than one and its right algebra is an Abelian algebra of dimension one.   We further show that all the solvable algebras are isomorphic if their left nilpotent algebras are Heisenberg algebras of fixed dimension and their right algebras are Abelian algebras of dimension one.


\end{abstract}


\section{Introduction}

 Lie algebra is a branch of algebra, closely related with the Lie group theory, has various applications in physics and mechanics \cite{Ia15, SW14}.
Sophus Lie (1842-1899) started the research of Lie group and Lie algebra \cite{LS93}.  Later, Friedrich Engel (1861-1941), Wilhelm Killing (1847-1923), \'{E}lie Cartan (1869-1951), Hermann Weyl (1885-1955), Bartel van der Waerden (1903-1996), Nathan Jacobson (1910-1999) and Eugene Dynkin (1924-2014) made various contributions to the developments of Lie algebra \cite{EW06, GK96, Ia15, Ja79, SW14}.    In particular, Cartan \cite{Ca94} characterized solvable Lie algebras and semisimple algebras by their Killing forms and  completely classified semisimple Lie algebras.
Eugenio Elia Levi \cite{Le05} showed that a finite-dimensional Lie algebra is a semidirect sum of a solvable ideal and a semisimple  Lie subalgebra.  As semisimple Lie algebras were completely classified, people turn their attentions from the classification of general finite-dimensional Lie algebras to the classification of solvable Lie algebras and its subclass nilpotent Lie algebras \cite{GK96, RWZ88, SW14}.  The problem is completely solved up to dimension $7$, and the difficulty increases rapidly as the dimension increases \cite{SW14}.

The 1996 book of Goze and Khakimdjanov \cite{GK96} devoted to the study of nilpotent Lie algebras.   Some classes of nilpotent Lie algebras, such as Abelian Lie algebras, Heisenberg algebras, two-step algebras, filiform Lie algebras, characteristically nilpotent Lie algebras, nilpotent standard algebras, were studied in \cite{GK96}.

How can solvable nonnilpotent Lie algebras be classified?   Starting from 1993 \cite{RW93}, Winternitz and his collaborators  pursued such an approach:  to take a class of nilpotent Lie algebras and construct all extensions of these algebras to solvable ones \cite{NW94, SW09, SW12}.  Also see \cite{WLD08}. These were further put together in the  2014 book of \u{S}nobl and Winternitz \cite[Part 3]{SW14}.   For a finite dimensional Lie algebra, both the largest solvable ideal and the largest nilpotent ideal exist.  They are called the radical and the nilradical of the Lie algebra.    The Levi decomposition states that every finite dimensional Lie algebra is isomorphic to a semidirect sum of its radical and a semisimple Lie algebra.   The approach of \u{S}nobl and Winternitz \cite{SW14} uses the nilradical of a solvable nonnilpotent Lie algebra as the starting point to study such a solvable Lie algebra.   A rational algorithm for identifying the nilradical of a Lie algebra was presented in \cite[Chapter 7]{SW14}.

As discussed in \cite[Chapter 3]{EW06}, even in the three dimensional case, if the dimension of the derived algebra is two, then there are infinitely many non-isomorphic Lie algebras.  This example will be discussed as Example 2 in Section 3 of this paper.  They are solvable nonnilpotent Lie algebras.   Therefore, the classification work of solvable nonnilpotent Lie algebras can only divide such Lie algebras to classes.  Some classes may contain infinitely many non-isomorphic Lie algebras.

In this paper, we propose to decompose a finite dimensional solvable nonnilpotent Lie algebra to two nilpotent Lie algebras.   The starting point of such an approach is the smallest ideal of the lower central series of a solvable nonnilpotent Lie algebra.   In \cite{Qi19}, we studied three minimax ideal relations of a Lie algebra.   The smallest ideal of the lower central series is at one of the three minimax places of a Lie algebra.   We introduced near perfect ideas for a finite dimensional Lie algebra in \cite{Qi19}.   The largest near perfect ideal of a finite dimensional Lie algebra always exists and is exactly the  smallest ideal of the lower central series of that Lie algebra.   We call it the near perfect radical of that Lie algebra.   A Lie algebra is nilpotent if and only if its near perfect radical is zero.  The factor algebra of a Lie algebra by its near perfect radical is always nilpotent.   For a solvable Lie algebra, its near perfect radical is nilpotent.   Using these facts, we decompose a finite dimensional solvable nonnilpotent Lie algebra to two nilpotent Lie algebras.   The first nilpotent Lie algebra is the factor algebra of the solvable Lie algebra by its near perfect radical.  We call it the left nilpotent algebra of that solvable Lie algebra.   The second Lie algebra is the near perfect radical of that solvable Lie algebra.  We call it the right nilpotent algebra of that solvable Lie algebra.  Then we may classify solvable Lie algebras by such a nilpotent decomposition, and study solvable Lie algebras with different nilpotent decompositions separately.    For example, we may study A($n_1$)-A($n_2$) algebras, which are solvable Lie algebras whose left and right nilpotent algebras are both Abelian and with dimensions $n_1$ and $n_2$ respectively.   Actually, the simplest solvable nonnilpotent Lie algebra is the A($1$)-A($1$) algebra, which is of dimension $2$.    See Example 1 in Section 3.  On the other hand, the three dimensional case stated in the last paragraph is actually the A($1$)-A($2$) algebra, which will be studied as Example 2 in Section 3.

We also combine the lower central series of the left and right nilpotent algebras of a solvable nonnilpotent Lie algebra together, and call it the extended lower central series of the solvable Lie algebra.  This strengthens the lower central series of such a solvable Lie algebra.

In the next section, some preliminaries are given.

In Section 3, we introduce the nilpotent decomposition and the extended lower central series for a finite dimensional solvable nonnilpotent Lie algebra, and study such a decomposition for dimension up to $4$.

We study A($n_1$)-A($n_2$) algebras in Section 4.  We show that an $A(n)$-$A(1)$ algebra with $n \ge 2$ is the direct sum of an A$(1)$-A$(1)$ algebra and $n-1$ Abelian algebras A$(1)$, and conjecture that an A($n_1$)-A($n_2$) algebra with $n_1 > n_2$ is always decomposable.   Then we give a theorem for the structure of general A($n_1$)-A($n_2$) algebras, for further classifying such solvable Lie algebras.

In Section 5, we study A($n$)-H$_p$ algebras, H$_p$-A($n$) algebras and H$_p$-H$_q$ algebras. Here, H$_p$ is the $p$th order Heisenberg algebra, which is of dimension $2p+1$.  We show that all the H$_p$-A($1$) algebras are isomorphic.   We then give a theorem for the structure of general H$_p$-A($n$) algebras.

Thus, we completely classify A$(1)$-A$(1)$ algebras and H$_p$-A($1$) algebras, and give outlines for classifying A($n_1$)-A($n_2$) algebras and H$_p$-A($n$) algebras.


Some final remarks are made in Section 6.

\section{Preliminaries}
\label{sec:Lies}

 In this paper, we only study finite dimensional Lie algebras.  The related field $\F$ is either the field of complex numbers or the field of real numbers.   Suppose that $L$ is an $n$-dimensional Lie algebra defined on $\F$.

 Denote the Lie bracket operation on $L$  by $[ \cdot, \cdot ]$.  Then $[ \cdot, \cdot ]$ is a bilinear operation, and for $\x, \y ,\z \in L$, we have
 \begin{equation} \label{e2.1}
 [\x, \y] = - [\y, \x]
 \end{equation}
 and the Jacobi identity
 \begin{equation} \label{e2.2}
 [\x, [\y, \z]] + [\y, [\z, \x]] + [\z, [\x, \y]] = \0.
 \end{equation}

 If $[\x, \y] \equiv 0$ for any $\x, \y \in L$, then $L$ is called an Abelian algebra.  We denote an Abelian algebra of dimension $n$ by A$(n)$.   When $n > 1$, Abelian algebra A$(n)$ can always be decomposed to a direct sum of $n$ Abelian algebras A$(1)$.  In \cite[Part 4]{SW14}, A$(1)$ is denoted as $n_{1,1}$.

 A subspace $I$ of $L$
 is called a subalgebra of $L$ if for any $\x, \y$ in $I$, $[\x, \y] \in I$.   If furthermore for any $\x \in L$ and $\y \in I$, we have $[\x, \y] \in I$, then $I$ is called an ideal of $L$.   The center of $L$, defined as $Z(L) = \{ \x \in L : [\x, L] = 0 \}$, is an ideal of $L$.

 Suppose that $I$ and $J$ are two ideals of $L$.  Then $I \cap J$,
 $$I+J := \{ \x + \y : \x \in I, \y \in J \}$$
 and
 $$[I, J] := {\rm Span}\{ [\x, \y] : \x \in I, \y \in J \}$$
 are also ideals of $L$.   The derived algebra of $L$ is defined as $D(L) := [L, L]$.
 If $L = D(L)$, then $L$ is called a perfect Lie algebra.

 Suppose that $I$ is an ideal of $L$.  Then the quotient vector space $L/I = \{ \x + I : \x \in L \}$ is a Lie algebra with a Lie bracket on $L/I$ defined by
$$[\x + I, \y + I] := [\x, \y] + I, \ \ \forall \x, \y \in L,$$
and is called the quotient or factor algebra of $L$ by $I$.

 Let
 $$L^{(0)} = L, \ \ \ L^{(k+1)} = D(L^{(k)}).$$
 Then we have the derived series of $L$: $L^{(0)} \supseteq L^{(1)} \supseteq L^{(2)} \supseteq \cdots$.
 $L^{(k)}$ are ideals of $L$ for all $k$.   If for some $m$, $L^{(m)} = 0$, then $L$ is called a solvable Lie algebra.   An ideal $I$ of $L$ is called a solvable ideal of $L$ if it itself is a solvable Lie algebra.  If $I$ and $J$ are two solvable ideals of $L$, then $I+J$ is also a solvable ideal of $L$.   Since $0$ is a trivial solvable ideal of $L$, the largest solvable ideal of $L$ exists, is called the radical of $L$ and denoted as $R(L)$.   A nonzero Lie algebra $L$ is called a semisimple Lie algebra if it has no nonzero solvable ideals.  Then, a nonzero Lie algebra is semisimple if and only if its radical is zero.  A semisimple Lie algebra is always a perfect Lie algebra but not vice versa.

 Let
 $$L^0 = L, \ \ \ L^{k+1} = [L, L^k].$$
 Then we have the lower central series of $L$: $L^0 \supseteq L^1 \supseteq L^2 \supseteq \cdots$. $L^k$ are ideals of $L$ for all $k$.   If for some $r$, $L^r = 0$, then $L$ is called a nilpotent Lie algebra, and $r$ is called its nilindex.  A nilpotent Lie algebra is always a solvable Lie algebra but not vice versa. A Lie algebra $L$ is solvable if and only if $D(L)$ is nilpotent.  An ideal $I$ of $L$ is called a nilpotent ideal of $L$ if it itself is nilpotent.   If $I$ and $J$ are two nilpotent ideals of $L$, then $I+J$ is also a nilpotent ideal of $L$.   Since $0$ is a trivial nilpotent ideal of $L$, the largest nilpotent ideal of $L$ exists, and is called the nilradical of $L$ and denoted as $N(L)$.   Since a nilpotent
Lie algebra is solvable, we always have $N(L) \subseteq R(L)$.   By Theorem 13 of \cite{Ja79}, we have $[L, R(L)] \subseteq N(L)$.   Hence, we always have
$[L, R(L)] \subseteq N(L) \subseteq R(L)$.   For a solvable Lie algebra, the dimension of its nilradical is not less than a half of that solvable Lie algebra \cite[Page 99]{SW14}. This forms the base of studying all solvable extensions of a class of nilpotent Lie algebras \cite{NW94, RW93, SW09, SW12, SW14}.

A nilpotent Lie algebra $L$ with its nilindex equal to $2$ is called a two-step algebra.  A special two-step algebra is the Heisenberg algebra H$_p$ with dimension $2p+1$ \cite[Page 43]{GK96}, defined by a basis $\{ \x_1, \cdots, \x_{2p+1} \}$, with the Lie bracket operations:
$$[\x_1, \x_2] = [\x_3, \x_4] = \cdots = [\x_{2p-1}, \x_{2p}] = \x_{2p+1}.$$
Heisenberg algebras $H_p$ have wide applications in physics \cite{GK96, SW14}.  The Heisenberg algebra H$_1$ is denoted as $n_{3,1}$ in \cite[Part 4]{SW14}.  The Heisenberg algebra H$_2$ is denoted as $n_{5,3}$ in \cite[Part 4]{SW14}.

An $n$-dimensional nilpotent Lie algebra $L$ is called a filiform Lie algebra \cite[Page 40]{GK96} if the dimensions of the ideals in its lower central series obey dim$(L^k) = n-k-1$ for $1 \le k \le n-1$.  We denote such a filiform algebra as F$(n)$.   The smallest dimension of Filiform algebras is $4$.  In \cite[Part 4]{SW14}, F$(4)$ is denoted as $n_{4, 1}$.   On the other hand, there are two filiform algebras of dimension $5$.   They are denoted as $n_{5,5}$ and $n_{5,6}$ in \cite[Part 4]{SW14}.  We may denote them as F$(5)_1$ and F$(5)_2$ respectively.
There are five filiform algebras of dimension $6$. They are denoted as $n_{6,18}$, $n_{6,19}$, $n_{6,20}$, $n_{6,21}$ and $n_{6,22}$ in \cite[Part 4]{SW14}.   We may denote them as F$(6)_1$, F$(6)_2$,  F$(6)_3$, F$(6)_4$ and F$(6)_5$ respectively.

In \cite{Qi19}, near perfect ideals were introduced.

An ideal $I$ of $L$ is called a near perfect ideal of $L$ if $I = [L, I]$.  Since zero is a trivial near perfect ideal, $L$ always has a near perfect ideal.  The sum $I+J$ of two near perfect ideals $I$ and $J$ of $L$ is still a near perfect ideal of $L$.   Thus, the largest near perfect ideal of $L$ exists, and is called the near perfect radical of $L$, and denoted as $NP(L)$.  The near perfect radical $NP(L)$, is exactly the smallest ideal of the lower central series of $L$.
The factor algebra $L/NP(L)$ is always nilpotent. We always have $P(L) \subseteq NP(L)$.   A nonzero Lie algebra is nilpotent if and only if its near perfect radical is zero.

 \section{The Extended Lower Central Series}
\label{sec:ideal}

 Suppose that a finite dimensional Lie algebra $L$ is not nilpotent.  Then its lower central series has the form:
 \begin{equation} \label{e3.1}
L = L^0 \supsetneq L^1 \supsetneq \cdots \supsetneq L^{r} \not = 0,
\end{equation}
where $L^r = NP(L) \not = 0$.

\begin{proposition} \label{p3.1}
Suppose a finite dimensional Lie algebra $L$ is not nilpotent, and its lower central series has the form (\ref{e3.1}).  Then the factor algebra $L/NP(L)$ is nilpotent and the lower central series of $L/NP(L)$ has the form:
\begin{equation} \label{e3.2}
L/NP(L) \supsetneq L^1/NP(L) \supsetneq \cdots \supsetneq L^{r}/NP(L) = 0.
\end{equation}
If furthermore $L$ is solvable, then $NP(L)$ is also nilpotent.
\end{proposition}
{\bf Proof} By \cite[Proposition 3.5]{Qi19}, the factor algebra $L/NP(L)$ is nilpotent. Then (\ref{e3.2}) follows the correspondence between ideals of a Lie algebra \cite[Page 14]{EW06}.  If furthermore $L$ is solvable, then its derived algebra $D(L)$ is nilpotent.   Hence, its near perfect radical $NP(L)$ is also nilpotent, as $NP(L)$ is an ideal of $D(L)$. \qed

Suppose that $L$ is solvable nonnilpotent.   Assume the lower central series of $NP(L)$ has the form
 \begin{equation} \label{e3.3}
NP(L) = L^{, 0} \supsetneq L^{, 1} \supsetneq \cdots \supsetneq L^{, r_1} = 0.
\end{equation}

Then we have the extended lower central series of $L$
\begin{equation} \label{e3.4}
L = L^0 \supsetneq L^1 \supsetneq \cdots \supsetneq L^{r} = NP(L) = L^{, 0} \supsetneq L^{, 1} \supsetneq \cdots \supsetneq L^{, r_1} = 0.
\end{equation}
We call $L/NP(L)$ the left nilpotent algebra of $L$, and $NP(L)$ the right nilpotent algebra of $L$, and denote them as $LN(L)$ and $RN(L)$ respectively.   If $LN(L) \equiv L/NP(L)$ is an Abelian algebra A$(n_1)$, and $RN(L) \equiv NP(L)$ is an Abelian algebra A$(n_2)$, then we say that the solvable Lie algebra $L$ is an A$(n_1)$-A$(n_2)$ algebra.   Similarly, we have A$(n_1)$-H$_p$ algebras, H$_p$-A$(n_1)$ algebras, H$_{p_1}$-H$_{p_2}$ algebras, A$(n_1)$-F$(n_2)$ algebras, etc.

\medskip

In the examples, we always have $[\x_i, \x_j] = - [\x_j, \x_i]$ and $[\x_i, \x_j] = 0$ if both $[\x_i, \x_j]$ and $[\x_j, \x_i]$ are not defined.   For simplicity, from now on we discuss Lie algebras on the field of complex numbers only.

\medskip

{\bf Example 1} This is $s_{2, 1}$ in \cite[Chapter 16]{SW14}, the simplest solvable nonnilpotent Lie algebra.  The two dimensional solvable Lie algebra $L$ has a basis $\{ \x_1, \x_2 \}$.   We have $[\x_2, \x_1] = \x_1$.
Then we find the extended lower central series of $L$ has the form
$$L \equiv L^0 \supsetneq L^1 \equiv NP(L) \equiv L^{,0} \supseteq L^{,1} \equiv 0,$$
where
$$L^1 \equiv NP(L) \equiv L^{,0} = {\rm Span}\{ \x_1 \}.$$
According to our notation, this is an A$(1)$-A$(1)$ algebra.  In dimension $2$, this is the only solvable nonnilpotent Lie algebra.

\medskip

{\bf Example 2} These are $s_{3, 1}$ and $s_{3, 2}$ in \cite[Chapter 16]{SW14}.  The three dimensional solvable Lie algebra $L$ has a basis $\{ \x_1, \x_2, \x_3 \}$.   We have $[\x_3, \x_1] = \x_1$.   For $s_{3, 1}$, we have $[\x_3, \x_2] = a\x_2$, where the values of the parameter $a$ satisfy $0 < |a| \le 1$ and if $|a| = 1$, then arg$(a) \le \pi$.
For $s_{3, 2}$, we have $[\x_3, \x_2] = \x_1 + \x_2$.
Then we find the extended lower central series of $L$ has the form
$$L \equiv L^0 \supsetneq L^1 \equiv NP(L) \equiv L^{,0} \supseteq L^{,1} \equiv 0,$$
where
$$L^1 \equiv NP(L) = N(L) \equiv L^{,0} = {\rm Span}\{ \x_1, \x_2 \}.$$
According to our notation, this is an A$(1)$-A$(2)$ algebra.    As in \cite[Chapter 16]{SW14}, since the Casimir invariants of $s_{3, 1}$ and $s_{3, 2}$ are different, $s_{3, 1}$ and $s_{3, 2}$ are not isomorphic. In fact, as discussed in \cite[Chapter 3]{EW06}, $s_{3, 1}$ itself contains infinitely many non-isomorphic algebras, with the parameter $a$ as a parameter.   Note that the nilpotent radical $N(L)$ is the same for $s_{3, 1}$ and $s_{3, 2}$.  It is A$(2)$.    There are no other three-dimensional solvable nonnilpotent Lie algebras in the complex field.

\medskip

{\bf Example 3} This is $s_{4, 1}$ in \cite[Chapter 17]{SW14}.  The four dimensional solvable Lie algebra $L$ has a basis $\{ \x_1, \x_2, \x_3, \x_4 \}$.   We have $[\x_4, \x_2] = \x_1$ and  $[\x_4, \x_3] = \x_3$.
Then we find the extended lower central series of $L$ has the form
$$L \equiv L^0 \supsetneq L^1 \supsetneq L^2 \equiv NP(L) \equiv L^{,0} \supseteq L^{,1} \equiv 0,$$
where
$$L^1 = {\rm Span}\{ \x_1, \x_3 \},$$
$$L^2 \equiv NP(L) \equiv L^{,0} = {\rm Span}\{ \x_3 \}.$$
According to our notation, this is an H$_1$-A$(1)$ algebra.

\medskip

{\bf Example 4} These are $s_{4, 2}$, $s_{4, 3}$ and $s_{4, 4}$ in \cite[Chapter 17]{SW14}.  The four dimensional solvable Lie algebra $L$ has a basis $\{ \x_1, \x_2, \x_3, \x_4 \}$.   We have $[\x_4, \x_1] = \x_1$.   For $s_{4, 2}$, we have $[\x_4, \x_2] = \x_1 + \x_2$ and $[\x_4, \x_3] = \x_2 + \x_3$.   For $s_{4, 3}$, we have $[\x_4, \x_2] = a\x_2$ and $[\x_4, \x_3] = b\x_3$,  where the values of the parameters $a$ and $b$ satisfy $0 < |b| \le |a| \le 1$ and if $(a, b) \not = (-1, -1)$.  If one or both equalities hold, there are some further restrictions \cite[Chapter 17]{SW14}.
For $s_{4, 4}$, we have $[\x_4, \x_2] = \x_1 + \x_2$ and $[\x_4, \x_3] = a\x_3$, where $a \not = 0$.
Then we find the extended lower central series of $L$ has the form
$$L \equiv L^0 \supsetneq L^1 \equiv NP(L) \equiv L^{,0} \supseteq L^{,1} \equiv 0,$$
where
$$L^1 \equiv NP(L) = N(L) \equiv L^{,0} = {\rm Span}\{ \x_1, \x_2, \x_3 \}.$$
According to our notation, these are A$(1)$-A$(3)$ algebras.    As in \cite[Chapter 16]{SW14}, since the Casimir invariants of $s_{4, 2}$, $s_{4, 3}$ and $s_{4, 4}$ are different, they are not isomorphic, and are divided to three classes.   Note that the nilpotent radical $N(L)$ is the same for $s_{4, 1}$, $s_{4, 2}$, $s_{4, 3}$ and $s_{4, 4}$.  It is A$(3)$.   However, $s_{4, 1}$ is an H$_1$-A$(1)$ algebra in Example 3, while $s_{4, 2}$, $s_{4, 3}$ and $s_{4, 4}$ are A$(1)$-A$(3)$ algebras.

\medskip

{\bf Example 5} These are $s_{4, 6}$, $s_{4, 8}$ and $s_{4, 10}$ in \cite[Chapter 17]{SW14}.  The four dimensional solvable Lie algebra $L$ has a basis $\{ \x_1, \x_2, \x_3, \x_4 \}$.   We have $[\x_2, \x_3] = \x_1$.   For $s_{4, 6}$, we have $[\x_4, \x_2] = \x_2$ and $[\x_4, \x_3] = - \x_3$.   For $s_{4, 8}$, we have $[\x_4, \x_1] = (1+a)\x_1$, $[\x_4, \x_2] = \x_2$ and $[\x_4, \x_3] = a\x_3$,  where the values of the parameter $a$ satisfy $0 < |a| \le 1$ if $|a| = 1$, then arg$(a) \le \pi$.
For $s_{4, 10}$, we have $[\x_4, \x_1] = 2\x_1$, $[\x_4, \x_2] = \x_2$ and $[\x_4, \x_3] = \x_2 + \x_3$.
Then we find the extended lower central series of $L$ has the form
$$L \equiv L^0 \supsetneq L^1 \equiv NP(L) \equiv L^{,0} \supsetneq L^{,1} \supsetneq L^{,2} \equiv 0,$$
where
$$L^1 = {\rm Span}\{ \x_1, \x_2, \x_3 \} \equiv NP(L) \equiv L^{,0},$$
$$L^{,1} = {\rm Span}\{ \x_1 \}.$$
According to our notation, these are A$(1)$-H$_1$ algebras.    As in \cite[Chapter 16]{SW14}, since the Casimir invariants of $s_{4, 6}$, $s_{4, 8}$ and $s_{4, 10}$ are different, they are not isomorphic, and are divided to three classes.   The nilpotent radical $N(L)$ is H$_1$ for $s_{4, 6}$, $s_{4, 8}$ and $s_{4, 10}$.

\medskip

{\bf Example 6} This is $s_{4, 11}$ in \cite[Chapter 17]{SW14}.  The four dimensional solvable Lie algebra $L$ has a basis $\{ \x_1, \x_2, \x_3, \x_4 \}$.   We have $[\x_2, \x_3] = \x_1$,    $[\x_4, \x_1] = \x_1$ and $[\x_4, \x_2] = \x_2$.
Then we find the extended lower central series of $L$ has the form
$$L \equiv L^0 \supsetneq L^1 \equiv NP(L) \equiv L^{,0} \supsetneq L^{,1} \equiv 0,$$
where
$$L^1 \equiv NP(L) = N(L) \equiv L^{,0} = {\rm Span}\{ \x_1, \x_2 \}.$$
According to our notation, these are A$(2)$-A$(2)$ algebras.
The nilpotent radical $N(L)$ is H$_1$, the same as in the cases of $s_{4, 6}$, $s_{4, 8}$ and $s_{4, 10}$ in Example 5.   This shows that our classification is somewhat finer in this case.

\medskip

 There are no other four-dimensional non-decomposable, solvable nonnilpotent Lie algebras in the field of complex numbers.

\section{A($n_1$)-A($n_2$) Algebras}

As we discussed in Example 1, the A($1$)-A($1$) Algebra, i.e., $s_{2, 1}$ in \cite[Chapter 16]{SW14}, is the simplest solvable nonnilpotent Lie algebra.   In Examples 2, 4 and 6, we have studied the A$(1)$-A$(2)$ algebra, the A$(1)$-A$(3)$ algebra and the A$(2)$-A$(2)$ algebra, respectively.    In this section, we study more general A($n_1$)-A($n_2$) algebras.

Note that if $LN(L) = $ A$(n_1)$, then $RN(L) = D(L)$.

\subsection{A($1$)-A($n$) Algebras}

In this case, $L$ is an $(n+1)$-dimensional solvable nonnilpotent Lie algebra.  Its left nilpotent algebra $LN(L)$ is a one-dimensional Abelian algebra, and its right nilpotent algebra $RN(L)$ is an $n$-dimensional Abelian algebra.  Furthermore, $N(L) = D(L) = RN(L) =$ A$(n)$ in this case.  Hence, this falls to the case of solvable Lie algebras with one nonnilpotent element and an $n$-dimensional Abelian nilradical, and has been studied in \cite[Section 10.4]{SW14}.    For dimension $5$, $s_{5,5}, s_{5,6}, s_{5,7}, s_{5,9}, s_{5,10}$ in \cite[Chapter 18]{SW14} are A($1$)-A($4$) algebras.   For dimension $6$, $s_{6,10}, s_{6,11}, s_{6,12}, s_{6,13}, s_{6,14}, s_{6,17}, s_{6,18}$ in \cite[Chapter 19]{SW14} are A($1$)-A($5$) algebras.

\subsection{A($n$)-A($1$) Algebras}

An A($n$)-A($1$) Algebra has dimension $n+1$.   We now have the following theorem.

\begin{theorem} \label{t4.1}
Suppose that we have an $A(n)$-$A(1)$ algebra $L$, and $n \ge 2$.   Then $L$ is the direct sum of an A$(1)$-A$(1)$ algebra and $n-1$ Abelian algebras A$(1)$.
\end{theorem}
{\bf Proof} We have $RN(L) = D(L) =$ A($1$).  Let $\y$ be a nonzero element of $RN(L) = D(L) =$ A($1$).    Since $D(L) =$ Span$\{ \y \}$, there must exist $\x \in L$ such that
$[\x, \y] = \alpha \y$ and $\alpha \not = 0$.   By scaling, we may assume that $[\x, \y] = \y$.   Let $I =$ Span $\{ \x, \y \}$.  Then $I$ is an A($1$)-A($1$) algebra.   Since $n \ge 2$, the dimension of $L$ is $n+1 \ge 3$.   Hence, there exist $\uu_1, \cdots, \uu_{n-1} \in L$ such that $\{ \x, \y, \uu_1, \cdots, \uu_{n-1} \}$ are linearly independent.   Since $D(L) =$ Span$\{ \y \}$, we have
$$[\uu_j, \y] = \alpha_j \y, \ {\rm and}\ [\uu_j, \x] = \beta_j \y,$$
for $j = 1, \cdots, n-1$.
Let $\z_j = \uu_j - \alpha_j \x + \beta_j \y$, for $j = 1, \cdots, n-1$.  Then
\begin{equation} \label{e4.7}
[\z_j, \y] = [\z_j, \x] = 0,
\end{equation}
for $j = 1, \cdots, n-1$.
Furthermore, $\{ \x, \y, \z_1, \cdots, \z_{n-1} \}$ are linearly independent.
Let $i, j = 1, \cdots, n-1$.
By the Jacobi identity (\ref{e2.2}), we have
$$[\x, [\z_i, \z_j]] + [\z_i, [\z_j, \x]] + [\z_j, [\x, \z_i]] = 0.$$
By (\ref{e4.7}), we have
\begin{equation} \label{e4.8}
[\x, [\z_i, \z_j]] = 0.
\end{equation}
Since $D(L) =$ Span$\{ \y \}$, we have
\begin{equation} \label{e4.9}
[\z_i, \z_j] = \gamma_{ij}\y.
\end{equation}
By (\ref{e4.8}) and (\ref{e4.9}), we have
$$0 =  [\x, [\z_i, \z_j]] = \gamma_{ij}[\x, \y] = \gamma_{ij}\y.$$
Hence,
\begin{equation} \label{e4.10}
[\z_i, \z_j] = 0
\end{equation}
for $i, j = 1, \cdots, n-1$.
By (\ref{e4.7}) and (\ref{e4.10}), $L$ is the direct sum of $I$ and Span$\{ \z_i\}$ for $i = 1, \cdots, n-1$.   The proof is completed.
\qed

\subsection{A($n_1$)-A($n_2$) Algebras for $n_2 \ge 2$}

In Examples 1, 2, 4 and 6, we discussed A($n_1$)-A($n_2$) algebras in \cite[Chapters 16-17]{SW14} for $n=n_1+n_2 \le 4$. In Subsection 4.1, we listed A($1$)-A($n$) algebras in \cite[Chapters 18-19]{SW14} for $n+1 = 5$ and $6$.

There are other A($n_1$)-A($n_2$) algebras in \cite[Chapters 18-19]{SW14} for $n = n_1+n_2 = 5$ and $6$.   For dimension $5$, $s_{5,41}$ and $s_{5,42}$ in \cite[Chapter 18]{SW14} are A($2$)-A($3$) algebras.   For dimension $6$, in \cite[Chapter 19]{SW14}, $s_{6,82}$, $s_{6,83}$, $s_{6,84}$, $s_{6,85}$, $s_{6,86}$, $s_{6,87}$, $s_{6,88}$, $s_{6,89}$, $s_{6,116}$, $s_{6,140}$, $s_{6,141}$, $s_{6,143}$, $s_{6,144}$, $s_{6,146}$, $s_{6,194}$, $s_{6,195}$ and $s_{6,196}$ are A($2$)-A($4$) algebras, $s_{6,143}$, $s_{6,183}$ and $s_{6,189}$ are A($3$)-A($3$) algebras.

In fact, in \cite[Chapters 16-19]{SW14}, we have not found any A($n_1$)-A($n_2$) algebra with $n_1 > n_2$.    Thus, we conjecture that Theorem \ref{t4.1} may be further extended to this case, i.e., we conjecture that any A($n_1$)-A($n_2$) algebra with $n_1 > n_2$ is decomposable.

When $n_2 \ge 2$, we may take a way similar to \cite[Section 10.1]{SW14} to classify A($n_1$)-A($n_2$) algebras.  We have the following theorem to describe the structure of such an $A(n_1)$-$A(n_2)$ algebra.

\begin{theorem} \label{t4.2}
Suppose that we have an $A(n_1)$-$A(n_2)$ algebra $L$, and $n_2 \ge 2$.   Assume that $\{ \y_1, \cdots, \y_{n_2} \}$ is a basis of $D(L)$, and there are vectors $\x_1, \cdots, \x_{n_1} \in L$ such that $\{ \y_1, \cdots, \y_{n_2}, \x_1, \cdots, \x_{n_1} \}$ is a basis of $L$. Then we have the following conclusions:
\begin{equation} \label{e4.11}
[\y_i, \y_j] = 0,
\end{equation}
\begin{equation} \label{e4.12}
[\x_\alpha, \y_i] = D_\alpha \y_i,
\end{equation}
\begin{equation} \label{e4.13}
[\x_\alpha, \x_\beta] = \sum_{k=1}^{n_2}b_{\alpha\beta}^k\y_k,
\end{equation}
for $1 \le i < j \le n_2$, $1 \le \alpha < \beta \le n_1$. Here $D_\alpha$ are $n_2 \times n_2$ matrices, map $D(L)$ to $D(L)$, and satisfy
\begin{equation} \label{e4.14}
[D_\alpha, D_\beta] \equiv D_\alpha D_\beta - D_\beta D_\alpha = 0,
\end{equation}
for $\alpha, \beta = 1, \cdots, n_1$. Furthermore,
\begin{equation} \label{e4.15}
{\rm Span}\{ [\x_\alpha, \y_i] : \alpha = 1, \cdots, n_1, i = 1, \cdots, n_2 \} = D(L)
\end{equation}
For $n_1 \ge 3$, the matrices $D_\alpha$ and the constants $b_{\alpha\beta}^k$ satisfy
\begin{equation} \label{e4.16}
\sum_{k=1}^{n_2} \left[ b_{\alpha\beta}^k (D_\gamma)_k^j + b_{\beta\gamma}^k (D_\alpha)_k^j + b_{\gamma\alpha}^k (D_\beta)_k^j \right] = 0,
\end{equation}
for $1 \le \alpha < \beta < \gamma \le n_1, 1 \le j \le n_2$.
A classification of the $A(n_1)$-$A(n_2)$ algebra $L$ thus amounts to a classification of the matrices $D_\alpha$ and the constants $b_{\alpha\beta}^k$ under the transformations:

(1) Reselect the representative vectors of the factor algebra $LN(L)$:
\begin{equation} \label{e4.17}
\tilde \x_\alpha = \x_\alpha + \sum_{j=1}^{n_2} R_\alpha^j \y_j.
\end{equation}

(2) Change the basis vectors $\{ \y_1, \cdots, \y_{n_2} \}$ of $D(L)$ by
\begin{equation} \label{e4.18}
\tilde \y_i = \sum_{j=1}^{n_2} S_i^j \y_j,
\end{equation}
where $(S_i^j)$ is an $n_2 \times n_2$ nonsingular matrix.

(3) Change the complementary basis vectors $\{ \x_1, \cdots, \x_{n_1} \}$ by
\begin{equation} \label{e4.19}
\tilde \x_\alpha = \sum_{\beta=1}^{n_1} G_\alpha^\beta \x_\beta,
\end{equation}
where $(G_\alpha^\beta)$ is an $n_1 \times n_1$ nonsingular matrix.
\end{theorem}
{\bf Proof} By the Abelian properties of $LN(L)$ and $RN(L)$, we have (\ref{e4.11}-\ref{e4.13}).  By the Jacobi identity for $\x_\alpha, \x_\beta, \y_i$, we have (\ref{e4.14}).  By the near perfect property of $D(L)$, we have (\ref{e4.15}).   By the Jacobi identity for $\x_\alpha, \x_\beta, \x_\gamma$, we have (\ref{e4.16}).   By considering all possible choices of $\{ \y_1, \cdots, \y_{n_2}, \x_1, \cdots, \x_{n_1} \}$, we have (\ref{e4.17}-\ref{e4.19}). \qed

We see that Theorem \ref{t4.2} is similar to Theorem 10.1 of \cite{SW14}.  There are differences.   First, we do not require linear nilindependence of the matrices $D_\alpha$ here as $D(L)$ is not the nilradical in general.    Second, we need (\ref{e4.15}) as $D(L) = NP(L)$ in our case.

Then $\{ D_1, \cdots, D_{n_1} \}$ is a basis of an Abelian subalgebra of the general linear algebra of $n_1 \times n_1$ matrices.   We may take an approach similar to \cite{SW14} to classify all such Abelian subalgebras and transform the matrices $ D_1, \cdots, D_{n_1}$ to some canonical forms.   Then we may determine the structure constants $b_{\alpha\beta}^k$, and weed out decomposable obtained A($n_1$)-A($n_2$) algebras.  We do not go to the details.

Hence, we see that A($n_1$)-A($n_2$) algebras are somewhat close to solvable extensions of Abelian nilradicals.  However, they contain different solvable Lie algebras, as the nilradical of an A($n_1$)-A($n_2$) algebra may not be Abelian, and a solvable extension of an Abelian nilradical may not be an A($n_1$)-A($n_2$) algebra.   This, in a certain sense, may deepen our understanding to solvable nonnilpotent algebras.

\section{A($n$)-H$_p$ Algebras, H$_p$-A($n$) Algebras and H$_p$-H$_q$ Algebras}

In Examples 3 and 5, we studied the H$_1$-A$(1)$ algebra and the  A$(1)$-H$_1$ algebra respectively.    In this section, we study more general A($n$)-H$_p$ algebras and H$_p$-A($n$) algebras.

\subsection{A($1$)-H$_p$ Algebras}

In this case, the nilradical of an A($1$)-H$_p$ algebra is the Heisenberg algebra H$_p$.
Then this is studied in \cite[Chapter 11]{SW14}.   In Example 5, we studied the  A$(1)$-H$_1$ algebra.   In \cite[Chapter 19]{SW14}, $s_{6,162}$, $s_{6,163}$, $s_{6,168}$, $s_{6,169}$, $s_{6,170}$, $s_{6,171}$, $s_{6,178}$, $s_{6,179}$, $s_{6,181}$ and $s_{6,182}$ are A($1$)-H$_2$ algebras.

\subsection{H$_p$-A($1$) Algebras}

In Subsection 4.2, we showed that an A($n$)-A($1$) algebra is always decomposable.  The next theorem shows that an H$_p$-A($1$) algebra is very different.

\begin{theorem} \label{t5.1}
Suppose that $L$ is an H$_p$-A($1$) algebra where $p$ is a positive integer.  Then $L$ has a basis $\{ \x_1, \cdots, \x_{2p+1}, \y \}$ such that $NP(L) =$ Span $\{ \y \}$,
\begin{equation} \label{e5.11}
[\x_1, \x_2] = \cdots = [\x_{2p-1}, \x_{2p}] = \x_{2p+1},
\end{equation}
\begin{equation} \label{e5.12}
[\x_1, \y] = \y,
\end{equation}
\begin{equation} \label{e5.13}
[\x_i, \y] = 0,  \ {\rm for}\ i=2, \cdots, 2p+1,
\end{equation}
and
\begin{equation} \label{e5.14}
[\x_i, \x_j] = 0,  \ {\rm for}\ i, j = 1, \cdots, 2p+1, (i,j) \not = (2l-1, 2l), \ {\rm for}\ l = 1, \cdots, p.
\end{equation}
Hence, all the H$_p$-A($1$) algebras are isomorphic.
\end{theorem}
{\bf Proof} Since $RN(L) = NP(L) =$ A$(1)$, we may assume that there exists a nonzero element $\y$ such that $RN(L) = NP(L) =$ $I =$ Span$\{ \y \}$.  Since $LN(L) = L/I = $ H$_p$, there exists a basis $\{ \uu_1 + I, \cdots, \uu_{2p+1} + I \}$ of $L/I$ such that $\uu_1, \cdots, \uu_{2p+1} \in L$,
\begin{equation} \label{e5.15}
[\uu_{2l-1}, \uu_{2l}] - \uu_{2p+1} \in I, \ {\rm for}\ l = 1, \cdots, p,
\end{equation}
\begin{equation} \label{e5.16}
[\uu_i, \y] = \alpha_i \y,  \ {\rm for}\ i=1, \cdots, 2p+1,
\end{equation}
and
\begin{equation} \label{e5.17}
[\uu_i, \uu_j] \in I,  \ {\rm for}\ i, j = 1, \cdots, 2p+1, (i,j) \not = (2l-1, 2l), \ {\rm for}\ l = 1, \cdots, p.
\end{equation}
By $[\uu_1, \uu_2] - \uu_{2p+1} \in I$, we have
$$0 = \left[[\uu_1, \uu_2] - \uu_{2p+1}, \y\right] = [[\uu_1, \y], \uu_2] + [\uu_1, [\uu_2, \y]] - [\uu_{2p+1}, \y] = -\alpha_{2p+1}\y.$$
This implies that $\alpha_{2p+1}=0$.
As $I = NP(L)$, by the properties of a near perfect ideal, $[L, I] = I$.   Therefore,
at least for one $i, \alpha_i \not = 0$, where $i=1, \cdots, 2p$.   Then we may reorder $i = 1, \cdots, 2p$, and make some scalings if necessary, such that $\alpha_1 = 1$ and the other parts of (\ref{e5.15}-\ref{e5.17}) still hold.

Let $\vv_1 = \uu_1$, $\vv_2 = \uu_2 -\alpha_2\uu_1 - \sum_{l=2}^p (\alpha_{2l}\uu_{2l-1}-\alpha_{2l-1}\uu_{2l})$, $\vv_i =\uu_i -\alpha_i\uu_1$ for $i = 3, \cdots, 2p$, and $\vv_{2p+1} = \uu_{2p+1}$.  Then we have
\begin{equation} \label{e5.18}
[\vv_{2l-1}, \vv_{2l}] = \vv_{2p+1} +\beta_{2l-1,2l}\y, \ {\rm for}\ l = 1, \cdots, p,
\end{equation}
\begin{equation} \label{e5.19}
[\vv_1, \y] = \y,
\end{equation}
\begin{equation} \label{e5.20}
[\vv_i, \y] = 0,  \ {\rm for}\ i=2, \cdots, 2p+1,
\end{equation}
and
\begin{equation} \label{e5.21}
[\vv_i, \vv_j] = \beta_{i,j}\y,  \ {\rm for}\ i, j = 1, \cdots, 2p+1, (i,j) \not = (2l-1, 2l), \ {\rm for}\ l = 1, \cdots, p.
\end{equation}

Let $\x_1 = \vv_1$, $\x_2 = \vv_2 - \beta_{1,2}\y - \beta_{1, 2p+1}\y$ and $\x_i =\vv_i -\beta_{1,i}\y$ for $i = 3, \cdots, 2p+1$.  Then we have
\begin{equation} \label{e5.22}
[\x_1, \x_2] = \x_{2p+1}, \ [\x_{2l-1}, \x_{2l}] = \x_{2p+1} +\sigma_{2l-1,2l}\y, \ {\rm for}\ l = 2, \cdots, p,
\end{equation}
where $\sigma_{2l-1,2l} = \beta_{2l-1,2l} + \beta_{1, 2p+1}$ for $l = 2, \cdots, p$,
\begin{equation} \label{e5.23}
[\x_1, \y] = \y,
\end{equation}
\begin{equation} \label{e5.24}
[\x_i, \y] = 0,  \ {\rm for}\ i=2, \cdots, 2p+1,
\end{equation}
\begin{equation} \label{e5.25}
[\x_1, \x_j] = 0, \ {\rm for}\ j = 3, \cdots,  2p+1,
\end{equation}
\begin{equation} \label{e5.26}
[\x_i, \x_j] = \sigma_{i,j}\y,  \ {\rm for}\ i, j = 2, \cdots, 2p+1, (i,j) \not = (2l-1, 2l), \ {\rm for}\ l = 2, \cdots, p,
\end{equation}
where $\sigma_{i,j} = \beta_{i,j}$ for $i, j = 2, \cdots, 2p+1, (i,j) \not = (2l-1, 2l)$ for $l = 2, \cdots, p$.

These are very close to the conclusions (\ref{e5.11}-\ref{e5.14}) except that we have to prove or make
$$\sigma_{i,j}=0,  \ {\rm for}\ i, j = 2, \cdots, 2p+1.$$
We now work on this.

Let $i, j = 3, \cdots, 2p+1$.   By the Jacobi identity and (\ref{e5.25}),
$$[\x_1, [\x_i, \x_j]] = [\x_i, [\x_1, \x_j]] - [\x_j, [\x_1, \x_i]] = 0.$$
On the other hand, by (\ref{e5.23}) and (\ref{e5.26}), we have
$$[\x_1, [\x_i, \x_j]] = \sigma_{i,j}\y.$$
These show that
\begin{equation} \label{e5.27}
\sigma_{i, j} = 0, \ {\rm for}\ i, j = 3, \cdots, 2p+1.
\end{equation}
Now, let $i, j = 2, \cdots, 2p+1$.   By the Jacobi identity and (\ref{e5.22}-\ref{e5.27}),
$$[\x_1, [\x_i, \x_j]] = [\x_i, [\x_1, \x_j]] - [\x_j, [\x_1, \x_i]] = 0.$$
On the other hand, by (\ref{e5.22}-\ref{e5.27}), we have
$$[\x_1, [\x_i, \x_j]] = \sigma_{i,j}\y.$$
These show that
$$\sigma_{i, j} = 0, \ {\rm for}\ i, j = 2, \cdots, 2p+1.$$
The proof is completed.
\qed

In Example 3, if we replace $\x_1, \x_3$ and $\x_4$ in the example by $\x_3, \y$ and $\x_1$, then we have the form of this theorem for $p=1$.

\medskip

{\bf Example 7} This is $s_{6, 26}$ in \cite[Chapter 19]{SW14}.  The six dimensional solvable Lie algebra $L$ has a basis $\{ \e_1, \e_2, \e_3, \e_4, \e_5, \e_6 \}$.   If we let $\x_1 = \e_6$, $\x_2 = \e_5$, $\x_3 = \e_2$, $\x_4 = \e_3$, $\x_5 = \e_1$ and $\y = \e_4$, then we have the form of this theorem for $p=2$.
Note that for an H$_2$-A($1$) algebra, the dimensions of the ideals in its lower central series should be CS$=[6, 2, 1]$.  In \cite[Chapter 19]{SW14}, there is only one solvable Lie algebra which has CS$=[6, 2, 1]$.  This confirms this theorem.    However, our theorem is also true for $p \ge 3$.   This also enlarges our knowledge about solvable nonnilpotent algebras.

\medskip

\subsection{H$_p$-A($n$) Algebras for $n \ge 2$}

We now discuss H$_p$-A($n$) algebras for $n \ge 2$.   In \cite[Chapter 18]{SW14}, $s_{5,20}$ and $s_{5,39}$ are H$_1$-A($2$) algebras.  Thus, we need to classify H$_p$-A($n$) algebras for $n \ge 2$.   The following is a theorem about the structure of H$_p$-A($n$) algebras for $n \ge 2$.

\begin{theorem} \label{t5.2}
Suppose that we have an H$_p$-A($n$) algebra $L$, and $n \ge 2$.   Assume that $I = NP(L) = RP(L)$ has a basis $\{ \y_1, \cdots, \y_n \}$, and $L/I$ has a basis $\{ \x_1+I, \cdots, \x_{2p+1}+I \}$.
Then we have the following conclusions:
\begin{equation} \label{e5.37}
[\y_i, \y_j] = 0
\end{equation}
for $i, j = 1, \cdots, n$,
\begin{equation} \label{e5.38}
\left([\x_\alpha, \y_1], [\x_\alpha, \y_2], \cdots, [\x_\alpha, \y_n]\right) = \left(\y_1, \y_2, \cdots, \y_n\right) D_\alpha
\end{equation}
for $\alpha =1, \cdots, 2p+1, i = 1, \cdots, n$, where $D_\alpha$ are $n \times n$ matrices, map $I$ to $I$, and satisfy
\begin{equation} \label{e5.39}
[D_\alpha, D_\beta] \equiv D_\alpha D_\beta - D_\beta D_\alpha = \left\{ \begin{array} {ll} 
D_{2p+1}, & {\rm if}\ (\alpha, \beta) = (2l-1, 2l), \ 1 \le l \le p,\\ 
0, & {\rm otherwise},
\end{array}\right.
\end{equation}
for $1 \le \alpha < \beta \le 2p+1$,
\begin{equation} \label{e5.40}
[\x_{2l-1}, \x_{2l}] = \x_{2p+1} + \sum_{k=1}^n b_{2l-1, 2l}^k\y_k
\end{equation}
for $l=1, \cdots, p$,
\begin{equation} \label{e5.41}
[\x_\alpha, \x_\beta] = \sum_{k=1}^n b_{\alpha\beta}^k\y_k,
\end{equation}
for $1 \le \alpha < \beta \le 2p+1, (\alpha, \beta) \not = (2l-1, 2l)$ for $l = 1, \cdots, p$. Furthermore,
\begin{equation} \label{e5.42}
{\rm Span}\{ [\x_\alpha, \y_i] : \alpha = 1, \cdots, 2p+1, i = 1, \cdots, n \} = I
\end{equation}
The matrices $D_\alpha$ and the constants $b_{\alpha\beta}^k$ satisfy
\begin{equation} \label{e5.43}
\sum_{k=1}^{n} \left[ b_{2l-1,2l}^k (D_\gamma)_k^j + b_{2l,\gamma}^k (D_{2l-1})_k^j - b_{2l-1, \gamma}^k (D_{2l})_k^j \right] + b_{\gamma, 2p+1}^j = 0, \ j = 1, \cdots, n,
\end{equation}
for $2l < \gamma \le 2p+1$ and $l = 1, \cdots, p$, where $b_{2p+1, 2p+1}^j = 0, \ j = 1, \cdots, n$;
\begin{equation} \label{e5.43.1}
\sum_{k=1}^{n} \left[ b_{\alpha,2l-1}^k (D_{2l})_k^j - b_{2l-1,2l}^k (D_\alpha)_k^j + b_{\alpha, 2l}^k (D_{2l-1})_k^j \right] - b_{\alpha, 2p+1}^j = 0, \ j = 1, \cdots, n,
\end{equation}
for $1 \le \alpha < 2l-1$ and $l = 1, \cdots, p$; 
\begin{equation} \label{e5.43.2}
\sum_{k=1}^{n} \left[ b_{\alpha\beta}^k (D_\gamma)_k^j + b_{\beta\gamma}^k (D_\alpha)_k^j - b_{\alpha\gamma}^k (D_\beta)_k^j \right] = 0, \ j = 1, \cdots, n,
\end{equation}
for $1 \le \alpha < \beta < \gamma \le 2p+1, (\alpha, \beta) \not = (2l-1, 2l), (\beta, \gamma) \not = (2l-1, 2l)$.
A classification of the H$_p$-A($n$) algebra $L$ thus amounts to a classification of the matrices $D_\alpha$ and the constants $b_{\alpha\beta}^k$ under the transformations:

(1) Reselect the representative vectors of the factor algebra $L/I$:
\begin{equation} \label{e5.44}
\tilde \x_\alpha = \x_\alpha + \sum_{j=1}^{n} R_\alpha^j \y_j.
\end{equation}

(2) Change the basis vectors $\{ \y_1, \cdots, \y_n \}$ of $I$ by
\begin{equation} \label{e5.45}
\tilde \y_i = \sum_{j=1}^n S_i^j \y_j,
\end{equation}
where $(S_i^j)$ is an $n \times n$ nonsingular matrix.

(3) Change the complementary basis vectors $\{ \x_1, \cdots, \x_{n_1} \}$ by
\begin{equation} \label{e5.46}
\tilde \x_\alpha = \sum_{\beta=1}^n G_\alpha^\beta \x_\beta,
\end{equation}
where $(G_\alpha^\beta)$ is an $2p+1 \times 2p+1$ nonsingular matrix.
\end{theorem}
{\bf Proof} By the Abelian properties of $I$ and the Heisenburg properties of $L/I$, we have (\ref{e5.37}), (\ref{e5.38}), (\ref{e5.40}) and (\ref{e5.41}).  By the Jacobi identity for $\x_\alpha, \x_\beta, \y_i$, we have (\ref{e5.39}).  By the near perfect property of $I$, we have (\ref{e5.42}).   By the Jacobi identity for $\x_\alpha, \x_\beta, \x_\gamma$, we have (\ref{e5.43}-\ref{e5.43.2}).   By considering all possible choices of $\{ \y_1, \cdots, \y_n, \x_1, \cdots, \x_{2p+1} \}$, we have (\ref{e5.44}-\ref{e5.46}). \qed

Similar to Theorem 4.2, $\{ D_1, \cdots, D_{2p+1} \}$ is a basis of an Abelian subalgebra of the general linear algebra of $n \times n$ matrices.  This property is true for any solvable Lie algebra whose right nilpotent algebra is Abelian.  We may take an approach similar to \cite{SW14} to classify all such Abelian subalgebras and transform the matrices $ D_1, \cdots, D_{2p+1}$ to some canonical forms.   Then we may determine the structure constants $b_{\alpha\beta}^k$, and weed out decomposable obtained H$_p$-A($n$) algebras.  We do not go to the details.

\medskip

\subsection{A($n$)-H$_p$ Algebras and H$_p$-H$_q$ Algebras}

In \cite[Chapter 18]{SW14}, $s_{5,15}$ is an A($2$)-H$_1$ algebra, $s_{6,24}$, $s_{6,30}$ and $s_{6,38}$ are H$_1$-H$_1$ algebras.   Further discussion on such algebras are needed.

\section{Final Remarks}

In this paper, we proposed a nilpotent decomposition classification approach for solvable nonnilpotent Lie algebras, as a complementary approach to the existing approach of solvable extension of a given nilpotent radical, as discussed in \cite{SW14}.  We may see that in some cases, such as demonstrated by Theorems \ref{t4.1}, \ref{t4.2}, \ref{t5.1} and \ref{t5.2}, this approach may reveal some structures of  solvable nonnilpotent Lie algebras, different from the the approach of solvable extension of a given nilpotent radical.
Hence, it is worth further exploring on this approach.

It is known that not all nilpotent Lie algebras can be nilradicals of some solvable Lie algebras \cite[Page 100]{SW14}.   In order to have nontrivial solvable extensions, a given nilpotent Lie algebra must possess at least one nonnilpotent derivation.   A nilpotent algebra which has only nilpotent derivations and consequently is not a nilradical of any solvable Lie algebra is called a characteristically nilpotent algebra.  A survey on characteristically nilpotent algebras can be found in \cite{AC01}.   It is wondered if a characteristically nilpotent algebra can be a left or right nilpotent algebra of a solvable nonnilpotent Lie algebra or not.  On the other hand, is there any condition on a given nilpotent Lie algebra to be the left or right nilpotent algebra of a solvable nonnilpotent Lie algebra?  These two issues may also worth further exploring.

\bigskip

{\bf Acknowledgment}   The author is thankful to Professors Chengming Bai, Huixiang Chen, Libor \u{S}nobl and Yisheng Song for their comments.


\begin{thebibliography}{99}

\bibitem{AC01} J.M. Ancochea and R. Campoamor, ``Characteristically nilpotent Lie algebras: A survey'', {\sl Extracta Mathematicae \bf 16} (2001) 153-210.

\bibitem{Ca94} E. Cartan, {\sl Sur la Structure des Groupes de Transformation Finis and Contnus}, Th\'{e}se, Paris (1894).

\bibitem{EW06} K. Erdamann and M.J. Wildon, {\sl Introduction to Lie Algebras},  Springer, New York (2006).

\bibitem{GK96} M. Goze and Y. Khakimdjanov, {\sl Nilpotent Lie Algebras}, Springer, New York (1996).

\bibitem{Ia15} F. Iachello, {\sl Lie Algebras and Applications}, Second Edition, Springer, New York (2015).

\bibitem{Ja79} N. Jacobson, {\sl Lie Algebras}, Dover Publication Inc., New York (1979).  Republication of the 1962 original.


\bibitem{Le05} E.E. Levi, ``Sulla structtura dei gruppi finiti e continui'', {\sl Atti. Accad. Sci. Torino \bf 40} (1905) 551-565.

\bibitem{LS93} S. Lie and G. Scheffers, {\sl Vorlesungen \'{u}ber Kontinuerliche Gruppen}, Leipzig (1893).

\bibitem{NW94} J.-C. Ndogmo and P. Winternitz, ``Solvable Lie algebras with abelian nilradicals'', {\sl J. Phys. \bf 27} (1994) 405-423.

\bibitem{Qi19} L. Qi, ``Three minimax ideal relations of Lie algebras'', arXiv:1901.10687v11, (2019).


\bibitem{RWZ88} D.W. Rand, P. Winternitz and H. Zassenhaus, ``On the identification of a Lie algebra given by its structure constants. I: Direct decompositions, Levi decompositions, and nilradicals'', {\sl Linear Algebra Appl. \bf 109} (1988) 197-246.

\bibitem{RW93} R.J. Rubin and P. Winternitz, ``Solvable Lie algebras with Heisenburg nilradicals'', {\sl J. Phys. \bf 26} (1993) 1123-1138.

\bibitem{SW09} L. \u{S}nobl and P. Winternitz, ``All solvable extensions of a class nilpotent Lie algebras of dimension $n$ and degree of nilpotency $n-1$'', {\sl J. Phys. \bf 42} (2009) 105201.

\bibitem{SW12} L. \u{S}nobl and P. Winternitz, ``Solvable Lie algebras with Borel nilradicals'', {\sl J. Phys. \bf 45} (2012) 095202.

\bibitem{SW14} L. \u{S}nobl and P. Winternitz, {\sl Classification and Identification of Lie Algebras}, AMS, Providence (2014).

\bibitem{WLD08} Y. Wang, J. Lin and S. Deng, ``Solvable Lie algebras with quasifiliform nilradicals'', {\sl Comm. Algebra \bf 36} (2008) 4052-4067.


























%

%

%

%

%

%

%

%

%

%

%

%

%

%

%

%

%

%

%

%

%

%

%

%

%

\end{thebibliography}
\end{document}